\documentclass[11pt,fleqn,twoside]{article}
\usepackage{amsfonts,amssymb,latexsym,graphics}
\makeatletter
\newcommand{\prava}[1]{\small\it
\begin{flushleft}
Copyright \copyright \ 2000 by  #1
\end{flushleft}}

\newcommand{\name}[1]{\begin{flushleft}
                       \LARGE \bf #1
                       \end{flushleft}\vspace{-3mm}}

\newcommand{\Author}[1]{\begin{flushleft}
                       \it #1 \end{flushleft}}

\newcommand{\Adress}[1]{\begin{flushleft}
                       \it #1 \end{flushleft}}

\newcommand{\Date}[1]{\begin{flushleft}
                      \small  \it #1 \end{flushleft}}

\newcommand{\ehkol}{Author \ name}
\newcommand{\ohkol}{Article \ name}
\renewcommand{\@evenhead}{
\hspace*{-3pt}\raisebox{-15pt}[\headheight][0pt]{\vbox{\hbox to \textwidth 
{\thepage \hfil \ehkol}\vskip4pt \hrule}}}
\renewcommand{\@oddhead}{
\hspace*{-3pt}\raisebox{-15pt}[\headheight][0pt]{\vbox{\hbox to \textwidth 
{\ohkol \hfil \thepage}\vskip4pt\hrule}}}
\renewcommand{\@evenfoot}{}
\renewcommand{\@oddfoot}{}

\newcommand{\be}{\begin{equation}}
\newcommand{\ee}{\end{equation}}
\newcommand{\ba}{\hspace*{-5pt}\begin{array}}
\newcommand{\ea}{\end{array}}

\newcommand{\ds}{\displaystyle}
\makeatother

\begin{document}
\thispagestyle{empty}
\setcounter{page}{1}

\renewcommand{\ehkol}{S.-J. Yu and K. Toda}
\renewcommand{\ohkol}{Lax Pairs, Painlev\'e Properties and Exact
Solutions}

\begin{flushleft}
\footnotesize \sf Journal of Nonlinear Mathematical Physics \qquad
2000, V.7, N~1, \pageref{yu-toda-fp}--\pageref{yu-toda-lp}.
\hfill {\sc Letter}
\end{flushleft}

\vspace{-5mm}

\renewcommand{\footnoterule}{}
{\renewcommand{\thefootnote}{}
 \footnote{\prava{S.-J. Yu and K. Toda}}}

\name{Lax Pairs, Painlev\'e Properties \\
and Exact Solutions of the Calogero\\
Korteweg-de Vries Equation \\
and a New {\mathversion{bold}$(2 + 1)$}-Dimensional Equation}\label{yu-toda-fp}

\Author{Song-Ju YU~$^\dag$ and Kouichi TODA~$^\ddag$}

\Adress{Department of Physics, Ritsumeikan University,
Kusatsu, Shiga, 525-7755, Japan\\
$\dag$~E-mail: fpc30017@se.ritsumei.ac.jp\\
$\ddag$~E-mail: sph20063@se.ritsumei.ac.jp}

\Date{Received May 4, 1999; Revised June 13, 1999; Accepted
July 14, 1999}

\begin{abstract}
\noindent
We prove the existence of a Lax pair for the
Calogero Korteweg-de Vries~(CKdV) 
equation.  Moreover, we modify the $T$ operator in the the Lax pair
of the 
CKdV equation, in the search of a $(2 + 1)$-dimensional case and
thereby propose a new 
equation in $(2 + 1)$ dimensions.  We named this the
$(2 + 1)$-dimensional CKdV 
equation.  We show that the CKdV equation as well as the
$(2+1)$-dimensional CKdV
equation are integrable in the sense that they possess the Painlev\'e 
property.  Some exact solutions are also constructed.
\end{abstract}

\renewcommand{\theequation}{\thesection.\arabic{equation}}
\setcounter{equation}{0}

\section{Introduction}

In this paper  we attempt to extend the Calogero Korteweg-de
Vries~(CKdV)
equation 
to a $(2 + 1)$-dimensional equation.  The CKdV equation is a 
$(1 + 1)$-dimensional nonlinear equation~\cite{w} of the form
\begin{equation}
\label{eq1}
w_t + \frac{1}{4} w_{xxx} + \frac{3w_x}{8w^2} +
\frac{3w_x^3}{8w^2} - \frac{3w_x w_{xx}}{4w} = 0. \label{ckdv}
\end{equation}
Pavlov constructed (\ref{eq1}) using the new method for the description 
of an inf\/inite set of dif\/ferential substitutions and the KdV 
modif\/ications~\cite{p}.  We brief\/ly describe how the CKdV equation was 
constructed by Pavlov.  The Lax pair of the KdV equation
\begin{equation}
u_t + \frac{1}{4} u_{xxx} + \frac{3}{2} u u_x = 0 \label{kdv}
\end{equation}
has the form
\be
L = \partial_x^2 + u, \label{kdvl} 
\ee
\be
T = \partial_x L + \frac{1}{2} u \partial_x - \frac{1}{4} u_x + \partial_t.
\label{kdvt}
\ee
\newpage

\noindent
Pavlov obtained an inf\/inite set of dif\/ferential substitutions and the KdV 
modif\/ications from the Taylor expansion of the linear system for~(\ref{kdvl}) 
and~(\ref{kdvt}) respectively (see~\cite{p}).  The f\/irst order of an inf\/inite set 
of dif\/ferential substitutions is the Miura transformation
\begin{equation}
u = v^2 + \sigma v_x, \qquad (\sigma = \pm i). \label{miura}
\end{equation}
After substitution of the Miura transformation~(\ref{miura}) into the f\/irst 
order KdV mo\-di\-f\/i\-ca\-tions, we obtain the modif\/ied~KdV (mKdV) equation
\begin{equation}
v_t + \frac{1}{4} v_{xxx} + \frac{3}{2} v^2 v_x = 0. \label{mkdv}
\end{equation}
This equation admits the Lax representation
\be
L = \partial_x^2 + 2 \sigma v \partial_x, \label{mkdvl} 
\ee
\be
T = \partial_x L + \sigma v \partial^2 - \left(\frac{3}{2} v^2 + \frac{1}{2} 
\sigma v_x\right) \partial_x + \partial_t. \label{mkdvt}
\ee
The representation (\ref{mkdvl}), (\ref{mkdvt}) can be obtained
from the Lax pair of the KdV 
equation (\ref{kdvl}), (\ref{kdvt}) by the gauge transformation~\cite{b}.  In 
the second order, an inf\/inite set of dif\/ferential substitutions
and the KdV 
modif\/ications, lead to the Miura type transformation
\begin{equation}
v = - \frac{1}{2w} (1 + \sigma w_x)\label{miuratype}
\end{equation}
and the CKdV equation~(\ref{ckdv}).  Hamiltonian structures for the CKdV 
equation are discussed in~\cite{p}.

This paper is organized as follows.  In Section~2, we construct a Lax pair of 
the CKdV equation~(\ref{ckdv}) and propose a new equation in $(2 + 1)$ 
dimensions by the extension of the~$T$ operator for the CKdV equation.
We named it 
the $(2 + 1)$-dimensional CKdV equation.  Moreover, another dimensional 
extension is performed by changing the~$L$ operator~\mbox{\cite{zs,d,kd}}~as
follows:
\begin{equation}
L \mapsto L + \partial_y. \label{zakd}
\end{equation}
A $(2 + 1)$-dimensional equation obtained by the abovemethod is,
however, reduced to 
the KP equation.  In Section~3, the CKdV equation and the 
$(2 + 1)$-dimensional CKdV equation are proved to be integrable in the sense 
that they possess the Painlev\'e property.  The solutions to these equations 
are constructed by the Miura transformation in Section~4.  Section~5 is 
devoted to discussions.

\renewcommand{\theequation}{\thesection.\arabic{equation}}
\setcounter{equation}{0}

\section{The Lax pairs of the CKdV equation \\
and the {\mathversion{bold}$(2 + 1)$}-dimensional CKdV equation}

We conjecture that a Lax pair of the CKdV equation~(\ref{ckdv}) is of 
the form
\be
L = \partial_x^2 + g[w] \partial_x + h[w], \label{conjckdvl} 
\ee
\be
T = \partial_x L + T' + \partial_t, \label{conjckdvt}
\ee
where $g[w]$, $h[w]$ are functions of $w$ and its $x$-derivatives, and $T'$ is 
an unknown operator.  We can f\/ix the form of $g[w]$, $h[w]$ and $T'$ by the 
condition that the Lax equation
\begin{equation}
[L, T] = 0 \label{laxeq}
\end{equation}
gives the CKdV equation.  The result is
\be
g[w] = \frac{\sigma}{w}, \label{gw} 
\ee
\be
h[w] = -\frac{1}{4w^2} - \frac{\sigma w_x}{2w^2}, 
\label{hw} 
\ee
\be
\ds T' = \frac{\sigma}{2w} \partial_x^2 - 
\frac{1}{2w^2} \partial_x - \frac{3\sigma}{16w^3} + 
\frac{w_x}{4w^3} 
 - \frac{\sigma w_x^2}{16 w^3} + \frac{\sigma w_{xx}}{8w^2}. 
\label{td}
\ee
Hence the Lax pair of the CKdV equation is expressed as
\be
 L = \partial_x^2 + \frac{\sigma}{w} \partial_x - 
\frac{1}{4w^2} - \frac{\sigma w_x}{2w^2}, \label{ckdvl} 
\ee
\be
 T = \partial_x L + \frac{\sigma}{2w} \partial_x^2 - 
\frac{1}{2w^2} \partial_x - \frac{3 \sigma }{16w^3} + 
\frac{w_x}{4w^3} - \frac{\sigma w_x^2}{16w^3}
 + \frac{\sigma w_{xx}}{8 w^2} + \partial_t. 
\label{ckdvt}
\ee

Next we construct a new equation in $(2 + 1)$ dimensions.  For that, we modify 
the above~$T$ operator to include another spatial dimension as follows
\begin{equation}
T = \partial_z L + T'' + \partial_t, \label{conj2p1ckdvt}
\end{equation}
where $L$ is the same $L$ operator~(\ref{ckdvl}) for the CKdV equation. The Lax
equation~(\ref{laxeq}) gives not only the form of $T''$ but also a new 
equation.  They are
\be
\ba{l}
\ds  T'' = \frac{1}{2} \sigma \partial_x^{-1} \left(\frac{1}{w}\right)_z 
\partial_x^2 - \frac{1}{2w} \partial_x^{-1} 
\left(\frac{1}{w}\right)_z \partial_x + \frac{\sigma w_{xz}}{8w^2}
- \frac{\sigma w_x w_z}{8w^3} 
\vspace{3mm}\\
\ds \qquad - \frac{\sigma}{8w^2} \partial_x^{-1} 
\left(\frac{1}{w}\right)_z
  + \frac{w_x}{4w^2} \partial_x^{-1} 
\left(\frac{1}{w}\right)_z - \frac{\sigma}{16w} \partial_x^{-1} 
\left(\frac{1}{w^2}\right)_z  + \frac{\sigma}{16w} \partial_x^{-1} 
\left(\frac{w_x^2}{w^2}\right)_z
\ea \label{tdd}
\ee
and
\be
\ba{l}
\ds  w_t + \frac{1}{4} w_{xxz} + \frac{w_z}{4w^2} + \frac{1}{8} w_x 
\partial_x^{-1} \left(\frac{1}{w^2}\right)_z 
\vspace{3mm}\\
\ds \qquad +  
\frac{w_x^2 w_z}{2w^2} - \frac{1}{8} w_x \partial_x^{-1} 
\left(\frac{w_x^2}{w^2}\right)_z 
 -  \frac{w_x w_{xz}}{2w} - \frac{w_{xx} w_z}{4w} = 0, 
\ea\label{2p1ckdv}
\ee
respectively.  We name the above  equation the $(2 + 1)$-dimensional CKdV 
equation.  It follows from~(\ref{conj2p1ckdvt}) and~(\ref{tdd}) that the Lax 
pair of $(2 + 1)$-dimensional CKdV equation is given by
\be
L = \partial_x^2 + \frac{\sigma}{w} \partial_x -  
\frac{1}{4w^2} - \frac{\sigma w_x}{2w^2}, \label{2p1ckdvl} 
\ee
\be
\ba{l}
\ds  T = \partial_z L + \frac{1}{2} \sigma \partial_x^{-1} 
\left(\frac{1}{w}\right)_z \partial_x^2 -  \frac{1}{2w} 
\partial_x^{-1} \left(\frac{1}{w}\right)_z \partial_x + \frac{\sigma w_{xz}}{8w^2} 
\vspace{3mm}\\
\ds \qquad - \frac{\sigma w_x w_z}{8w^3} - \frac{\sigma }{8w^2} 
\partial_x^{-1} \left(\frac{1}{w}\right)_z + \frac{w_x}{4w^2} 
\partial_x^{-1} \left(\frac{1}{w}\right)_z 
\vspace{3mm}\\
\ds \qquad - \frac{\sigma}{16w} \partial_x^{-1} 
\left(\frac{1}{w^2}\right)_z + \frac{\sigma}{16w} 
\partial_x^{-1} \left(\frac{w_x^2}{w^2}\right)_z + \partial_t. 
\ea \label{2p1ckdvt}
\ee
Equation (\ref{2p1ckdv}) and the Lax pair (\ref{2p1ckdvl}), (\ref{2p1ckdvt}) 
are reduced to the CKdV equation and the Lax pair of the CKdV equation in the 
case of $x = z$.  In \cite{fsty1,fty1,fty2}, we developed the construction 
method for higher-dimensional integrable equation. For example, we considered 
the Calogero--Bogoyavlenskij--Schif\/f~(CBS) equation~\cite{c1,c2,c3,cd,b1,b2,s1},
\begin{equation}
u_t + \frac{1}{4} u_{xxz} + u u_z + \frac{1}{2} u_x \partial_x^{-1} u_z = 0,
\label{cbs}
\end{equation}
and the modif\/ied Calogero--Bogoyavlenskij--Schif\/f (mCBS)~\cite{b2},
\begin{equation}
v_t + \frac{1}{4} v_{xxz} + v^2 v_z + \frac{1}{2} v_x \partial_x^{-1} \left(v^2\right)_z 
= 0. \label{mcbs}
\end{equation}
These equations admit the Lax representations, 
respectively~\cite{fty1,fty2},
\be
L = \partial_x^2 + u, \label{cbsl} 
\ee
\be 
T = \partial_z L + \frac{1}{2} \partial_x^{-1} u_z \partial_x - \frac{1}{4} 
u_z + \partial_t, \label{cbst}
\ee
and
\be
L = \partial_x^2 + 2 \sigma v \partial_x, \label{mcbsl} 
\ee
\be
\ba{l}
\ds T = \partial_z L + \sigma \partial_x^{-1} v_z \partial_x^2 
 + \left(\frac{1}{2} \partial_x^{-1} \left(v^2\right)_z - 2 v 
\partial_x^{-1} v_z - \frac{1}{2} \sigma v_z\right) \partial_x 
\vspace{3mm}\\
\ds \qquad  + \frac{\sigma}{4} v_{xz} + \frac{\sigma}{2} v \left(\partial_x^{-1} 
\left(v^2\right)_z \right) 
+ \sigma \partial_x^{-1} v_t + \partial_t, 
\ea\label{mcbst}
\ee
respectively.
We obtained the mCBS equation from the CBS equation using the same Miura 
transformation (\ref{miura}) that connects the KdV equation with the mKdV 
equation \cite{fsty1,fty1,fty2}.  We checked that the transformation 
(\ref{miuratype}) connects the mCBS equation (\ref{mcbs}) and the 
$(2 + 1)$-dimensional CKdV equation (\ref{2p1ckdv}), i.e.,
\be
\ba{l}
\ds  v_t + \frac{1}{4} v_{xxz} + v^2 v_z + \frac{1}{2} v_x \partial_x^{-1} 
\left(v^2\right)_z 
\vspace{3mm}\\
\ds \qquad = \left(\frac{1}{2w^2} (1 + \sigma w_x) - \frac{\sigma 
}{2w} \partial_x\right) \left\{w_t + \frac{1}{4} w_{xxz} +  
\frac{w_z}{4w^2} + \frac{1}{8} w_x \partial_x^{-1} 
\left(\frac{1}{w^2}\right)_z\right. 
\vspace{3mm}\\
\ds \qquad  \left.+  \frac{w_x^2 w_z}{4w^2} - \frac{1}{8} w_x \partial_x^{-1} 
\left(\frac{w_x^2}{w^2}\right)_z -  \frac{w_x w_{xz}}{2w} - 
\frac{w_{xx} w_z}{4w}\right\}. 
\ea\label{frommbsto2p1ckdv}
\ee
These results are depicted in Fig.~1.

\begin{figure}[t]
\begin{center}
\hspace*{-10mm}{\begin{picture}(100,80)
\put(-7,62){\makebox(24,10){(\ref{miura})}}
\put(113,62){\makebox(24,10){(\ref{miuratype})}}
\put(-7,11){\makebox(24,10){(\ref{miura})}}
\put(113,11){\makebox(24,10){(\ref{miuratype})}}
\put(-100,50){\framebox(90,15){KdV(\ref{kdv})}}
\put(-5,58){\vector(2,0){20}}
\put(20,50){\framebox(90,15){mKdV(\ref{mkdv})}}
\put(115,58){\vector(2,0){20}}
\put(140,50){\framebox(90,15){CKdV(\ref{ckdv})}}
\put(-55,45){\line(0,-1){5}}
\put(-55,38){\line(0,-1){5}}
\put(-55,31){\vector(0,-1){11}}
\put(65,45){\line(0,-1){5}}
\put(65,38){\line(0,-1){5}}
\put(65,31){\vector(0,-1){11}}
\put(185,45){\line(0,-1){5}}
\put(185,38){\line(0,-1){5}}
\put(185,31){\vector(0,-1){11}}
\put(-100,1){\framebox(90,15){CBS(\ref{cbs})}}
\put(-5,9){\vector(2,0){20}}
\put(20,1){\framebox(90,15){mCBS(\ref{mcbs})}}
\put(115,9){\vector(2,0){20}}
\put(140,1){\framebox(90,15){$(2 + 1)$ CKdV(\ref{2p1ckdv})}}
\end{picture}}
\end{center}

\noindent
{\bf Figure 1:} {\small The 
dimensional extensions are indicated by broken arrows.  
These broken arrows 
indicate the modif\/i\-ca\-tion of $T$ operators for the search of the 
$(2 + 1)$-dimensional case.  Full arrows mean the Miura transformations.  The 
mKdV equation~(6) and the mCBS equation~(25) are induced by the Miura 
transformation~(5) from the KdV equation~(2) and CBS equation~(24).  We construct 
exact solutions of the CKdV equation~(1) and the $(2 + 1)$-dimensional CKdV 
equation~(21) from the solution of the mKdV equation~(6) and the mCBS
equation~(25), by using
using the Miura type transformation~(9).}
\label{extensions}
\end{figure}

We also extend the CKdV equation via (\ref{zakd}) \cite{zs,d,kd}.  Namely we 
consider
\begin{equation}
L = \partial_x^2 + \frac{\sigma}{w} \partial_x -\frac{1}{4w^2} - 
\frac{\sigma w_x}{2w^2} + \partial_y. \label{z-2p1ckdvl}
\end{equation}
The $T$ operator corresponding to (\ref{z-2p1ckdvl}) should be of the form
\be
\ba{l}
\ds  T = \partial_x^3 + \frac{3 \sigma }{2w} \partial_x^2 + 
\left\{- \frac{3}{4w^2} - \frac{3\sigma w_x}{2w^2} - 
\frac{3}{4} \sigma \partial_x^{-1} \left(\frac{1}{w}\right)_y\right\} 
\partial_x  - \frac{\sigma}{8w^3} 
\vspace{3mm}\\
\ds \qquad  + \frac{3w_x}{4w^3}
+ \frac{3\sigma w_x^2}{4w^3} - \frac{3\sigma 
w_{xx}}{8w^2} + \frac{3 \sigma w_y}{8w^2} 
 + \frac{3}{8} \partial_y^{-1} \left\{\frac{1}{w} 
\partial_x^{-1} \left(\frac{1}{w}\right)_{yy}\right\}+ \partial_t. 
\ea\label{z-2p1ckdvt}
\ee
We can construct the following equation from the Lax equation with 
(\ref{z-2p1ckdvl}) and (\ref{z-2p1ckdvt}),
\be
\ba{l}
\ds w_t + \frac{1}{4} w_{xxx} +  \frac{3w_x^3}{2w^2} - 
\frac{3w_x w_{xx}}{2w} - \frac{3\sigma w_y}{4w} - \frac{3}{4} w^2 
\partial_x^{-1} \left(\frac{1}{w}\right)_{yy} 
\vspace{3mm}\\
\ds \qquad - \frac{3}{4} \sigma w_x \partial_x^{-1} \left(\frac{1}{w}\right)_y - 
\frac{3}{4} \sigma w^2 \left(\partial_y^{-1} \left\{\frac{1}{w} 
\partial_x^{-1} \left(\frac{1}{w}\right)_{yy}\right\}\right)_x = 0. 
\ea \label{z-2p1ckdv}
\ee
However, the above equation is reduced to the KP equation
\begin{equation}
\left(u_t + \frac{1}{4} u_{xxx} + \frac{3}{2} u u_x\right)_x + \frac{3}{4} u_{yy} = 0 
\label{kp}
\end{equation}
by the transformation
\begin{equation}
w = - \frac{\sigma}{2\partial_y^{-1} u_x}. \label{tftokp}
\end{equation}
It follows that we cannot construct a new $(2 + 1)$-dimensional
 equation by
this method.

In the previous papers \cite{fty1,fty2} we modif\/ied both $L$ and $T$
operators 
for the KdV equation in searching for a $(3 + 1)$-dimensional
equation.
However, 
the Lax equation was reduced to the $(2 + 1)$-dimensional equation.  This 
equation separated the f\/irst and second order equations for the KP 
hierarchy~\cite{hl}.  Let us apply the same procedure to the CKdV
equation and
search for a $(3 + 1)$-dimensional Lax pair.  That is, we consider the Lax 
pair~(\ref{z-2p1ckdvl}) and
\begin{equation}
T = \partial_z L + T''' + \partial_t. \label{conj3p1ckdvt}
\end{equation}
However, we cannot f\/ix the form of $T'''$ by the Lax equation and, therefore, 
cannot construct a new equation in $(3 + 1)$ dimensions from the Lax 
pair~(\ref{z-2p1ckdvl}) and (\ref{conj3p1ckdvt}).


\setcounter{equation}{0}
\section{Painlev\'e analysis for the CKdV equation\\
and the  {\mathversion{bold}$(2 + 1)$}-dimensional CKdV equation}

To prove the Painlev\'e property \cite{wtc,kjh}
of the CKdV equation~(\ref{ckdv}) 
and the $(2 + 1)$-dimensional CKdV equation~(\ref{2p1ckdv}), we rewrite these 
equations by the change of of variable
\begin{equation}
W = \frac{1}{w}, \label{Ww}
\end{equation}
so that
\be
 W^2 W_t + \frac{1}{4} W^2 W_{xxx} + \frac{3}{8} W_x^3 + \frac{3}{8} W^4 W_x - 
\frac{3}{4} W W_x W_{xx} = 0, \label{wckdv}
\end{equation}
\be
\ba{l}
\ds  W^3 W_x W_{xt} - W^3 W_{xx} W_t + \frac{1}{4} W^3 W_x W_{xxxz} - 
\frac{1}{4} W^3 W_{xx} W_{xxz}
 + \frac{3}{4} W W_x^2 W_{xx} W_z 
\vspace{3mm}\\
\ds \qquad  + \frac{3}{4} W W_x^3 W_{xz} - 
\frac{3}{4} W_x^4 W_z + \frac{3}{4} W^4 W_x^2 W_z 
 + \frac{1}{4} W^5 W_x W_{xz} - \frac{1}{4} W^5 W_{xx} W_z  
\vspace{3mm}\\
\ds \qquad 
 - \frac{1}{2} 
W^2 W_x^2 W_{xxz} - \frac{1}{4} W^2 W_x W_{xxx} W_z - \frac{1}{4} W^2 W_x W_{xx} W_{xz} + 
\frac{1}{4} W^2 W_{xx}^2 W_z = 0. 
\ea\hspace{-5.46pt}\label{w2p1ckdv}
\ee
The solutions to~(\ref{wckdv}) and~(\ref{w2p1ckdv}) have the
form
\begin{equation}
W \sim W_0 \gamma^{\alpha}. \label{leadingsol}
\end{equation}
Here $\gamma$ is single valued about an arbitrary movable singular
manifold and
$\alpha$ is a negative integer (leading order).  By using leading order 
analysis, we obtain
\begin{equation}
\alpha = - 1, \qquad W_0^2 + \gamma_x^2 = 0. \label{leading}
\end{equation}
Substituting
\begin{equation}
W = \sum_{j = 0} W_j \gamma^{j - 1} \label{expand}
\end{equation}
into (\ref{wckdv}) and (\ref{w2p1ckdv}), leads to the resonances 
of~(\ref{wckdv}), namely
\begin{equation}
j = -1, 1, 3, \label{resonancesckdv}
\end{equation}
and the resonances of~(\ref{w2p1ckdv}), namely
\begin{equation}
j = -1, 1, 2, 3. \label{resonances2p1ckdv}
\end{equation}
The resonance $j = - 1$ corresponds to the arbitrary singularity manifold 
$\gamma$.  We used $MATHEMATICA$~\cite{mathe} to handle the calculation for 
the existence of arbitrary functions at the above resonances (except for 
$j = - 1$).
We f\/ind that $W_1$, $W_3$ are arbitrary for equation (\ref{wckdv}), 
and $W_1$, $W_2$, $W_3$ are arbitrary for equation (\ref{w2p1ckdv}).
Thus the 
general solution $W$ to (\ref{wckdv}) and (\ref{w2p1ckdv}) admits 
a suf\/f\/icient number of arbitrary functions, thus satisfying the Painlev\'e 
property.  Therefore the CKdV equation and the $(2 + 1)$-dimensional CKdV 
equation are integrable.

\setcounter{equation}{0}
\section{Exact solutions to the CKdV equation\\
 and the {\mathversion{bold}$(2 + 1)$}-dimensional 
CKdV equation}

In the previous section, the integrability of the CKdV equation and the 
$(2 + 1)$-dimensional CKdV equation was shown by the use of the
Painlev\'e test. In 
this section we shall construct exact solutions of the CKdV equation and the 
$(2 + 1)$-dimensional CKdV equation.

The mKdV equation~(\ref{mkdv}) has the solutions~\cite{h}
\begin{equation}
v_N = \sigma \left\{\log\left(\frac{f_N}{g_N}\right)\right\}_x, 
\label{solofmkdv}
\end{equation}
where $f_N$ and $g_N$ can be expressed as
\be
f_N = 1 + \sum_{n=1}^N \sum_{{}_NC_n} \eta_{i_1 \cdots i_n} 
\exp(\lambda_{i_1} + \cdots + \lambda_{i_n}), \label{fn} 
\ee
\be
g_N = 1 + \sum_{n=1}^N \sum_{{}_NC_n} (-1)^n\eta_{i_1 \cdots i_n} 
\exp(\lambda_{i_1} + \cdots + \lambda_{i_n}), \label{gn} 
\ee
\be
\lambda_j = p_j x + r_j t + s_j, \hspace{0.5cm} r_j = - \frac{1}{4} p^3_j, 
\label{lambda} 
\ee
\be
\eta_{jk} = \frac{(p_j - p_k)^2}{(p_j + p_k)^2}, \label{etajk} 
\ee
\be
\eta_{i_1 i_2 \cdots i_{n-1} i_n} = \eta_{i_1,i_2} \cdots \eta_{i_1,i_n} 
\cdots \eta_{i_{n - 1},i_n}. \label{aijklm}
\ee
Here ${}_NC_n$ indicates summation over all possible combinations of $n$ 
elements taken from~$N$, and symbols $s_j$ always denote arbitrary constants.  
We can solve the Miura type transformation~(\ref{miuratype}) for $w$ using 
solutions to the mKdV equation~(\ref{solofmkdv}). The solutions to the CKdV 
equation~(\ref{ckdv}) are
\begin{equation}
w_N = \sigma \left(\frac{f_N}{g_N}\right)^2 \int 
\left(\frac{f_N}{g_N}\right)^{-2} dx + c \left(\frac{f_N}{g_N}\right)^2, 
\label{solofckdv}
\end{equation}
where $c$ is an integration constant.  The above integral factor is rewritten 
as
\begin{equation}
\int \left(\frac{f_N}{g_N}\right)^{-2} dx = x + \frac{H_N}{f_N},
\label{integral}
\end{equation}
where
\begin{equation}
H_N = 4 \sum_{n=1}^N \frac{f_{N-1}(\hat{n})}{p_n} \label{hn}
\end{equation}
and
\be
\ba{l}
\ds f_{N-1}(\hat{j})= 1 + e^{\lambda_1} + \cdots + e^{\lambda_{j-1}} + 
e^{\lambda_{j+1}} + \cdots + e^{\lambda_{n}}
+ \eta_{12} e^{\lambda_1 + \lambda_2} + \cdots
\vspace{1mm}\\
\ds \qquad  + \eta_{1j-1} e^{\lambda_1
+ \lambda_{j-1}}  + \eta_{1j+1} e^{\lambda_1 + \lambda_{j+1}} + 
\cdots 
+ \eta_{1N} e^{\lambda_1 + \lambda_{N}} + \cdots 
\vspace{1mm}\\
\ds \qquad  + 
\eta_{j-1 j+1} e^{\lambda_{j-1} + \lambda_{j+1}} + \cdots + \eta_{j-1 N} 
e^{\lambda_{j-1} + \lambda_{N}} + \cdots + \eta_{j+1 j+2} e^{\lambda_{j+1} + 
\lambda_{j+2}} + \cdots 
\vspace{1mm}\\
\qquad + \eta_{j+1 N} e^{\lambda_{j+1} + \lambda_{N}} + \cdots + 
\eta_{N-1 N} e^{\lambda_{N-1} + \lambda_{N}} + \cdots 
\vspace{1mm}\\
\ds \qquad+\eta_{1 2 \cdots j-1 j+1 \cdots N-1 N} e^{\lambda_1 + \cdots 
+ \lambda_{j-1} + \lambda_{j+1}+ \cdots + \lambda_{N}},
\ea
\ee 
that is, $f_{N-1}$ is of the same structure as $f_N$, except for the
$j$ index.  
Equation~(\ref{integral}) is dif\/ferentiated with respect to $x$, i.e.,
\begin{equation}
H_{N,x} f_N - H_N f_{N,x} + f_N^2 - g_N^2= 0. \label{hn2}
\end{equation}
We checked equation~(\ref{hn2}) up to $N = 6$ by the use of
$MATHEMATICA$.  Fig.~2 
shows the solution~(\ref{solofckdv}) with $N = 1$, $p_1 = 1$, $s_1 = 2$ and 
$c = 0$.  In Fig.~3, we depict the case of $N = 2$, $p_1 = 1$, $s_1 = 2$, 
$p_2 = 0.5$, $s_2 = 5$ and $c = 0$.

We obtain solutions of the $(2 + 1)$-dimensional CKdV equation~(\ref{2p1ckdv}) 
using the identical procedure as with the construction of
solutions~(\ref{solofckdv}).  
Therefore, the form of the solutions are the same as~(\ref{solofckdv}).  
The dif\/ference between solutions of the 
$(2 + 1)$-dimensional CKdV equation and  the CKdV equation, is the 
dimensional extension of~(\ref{lambda}):
\begin{equation}
\lambda_j = p_j x + q_j z + r_j t + s_j, \qquad r_j = - \frac{1}{4} 
p_j^2 q_j. \label{dimexlambda}
\end{equation}
The propagation of the solution to the $(2 + 1)$-dimensional CKdV equation with 
$N = 1$, $p_1 = 1$, $q_1 = 3$, $s_1 = 2$ and $c = 0$ is shown in Fig.~4.  
Fig.~5 shows the solution with $N = 2$, $p_1 = 1$, $q_1 = 3$, $s_1 = 0$, 
$p_2 = 0.5$, $q_2 = - 3$, $s_2 = 0$ and $c = 0$.

\setcounter{equation}{0}
\section{Conclusions}

In this paper, we obtained the Lax pair of the CKdV equation and searched 
for the Lax pair of the higher dimensional CKdV equation using three methods.  The 
f\/irst method is to modify the~$T$ operator for the Lax pair of the
CKdV equation.
We then have obtained
the \mbox{$(2 + 1)$}-dimensional CKdV equation~(\ref{2p1ckdv}) 
and the Lax pair~(\ref{2p1ckdvl}) and~(\ref{2p1ckdvt}). The second method is 
to modify the~$L$ operator.  We constructed the Lax pair~(\ref{z-2p1ckdvl}), 
(\ref{z-2p1ckdvt}) and the equation~(\ref{z-2p1ckdv}).  
Equation~(\ref{z-2p1ckdv}) is, however, reduced to the KP equation by the 
transformation~(\ref{tftokp}).  In the last method, we unified
the f\/irst and 
second methods.  Using this method we can expect a new $(3 + 1)$-dimensional 
equation.  It, however, gives no consistent Lax equation, unlike
the f\/irst 
and second methods.  We also discussed the Painlev\'e property and exact 
solutions of equation~(\ref{ckdv}) and equation~(\ref{2p1ckdv}), which
proves
that the equations are integrable.


\newpage

\centerline{\scalebox{1}{\includegraphics{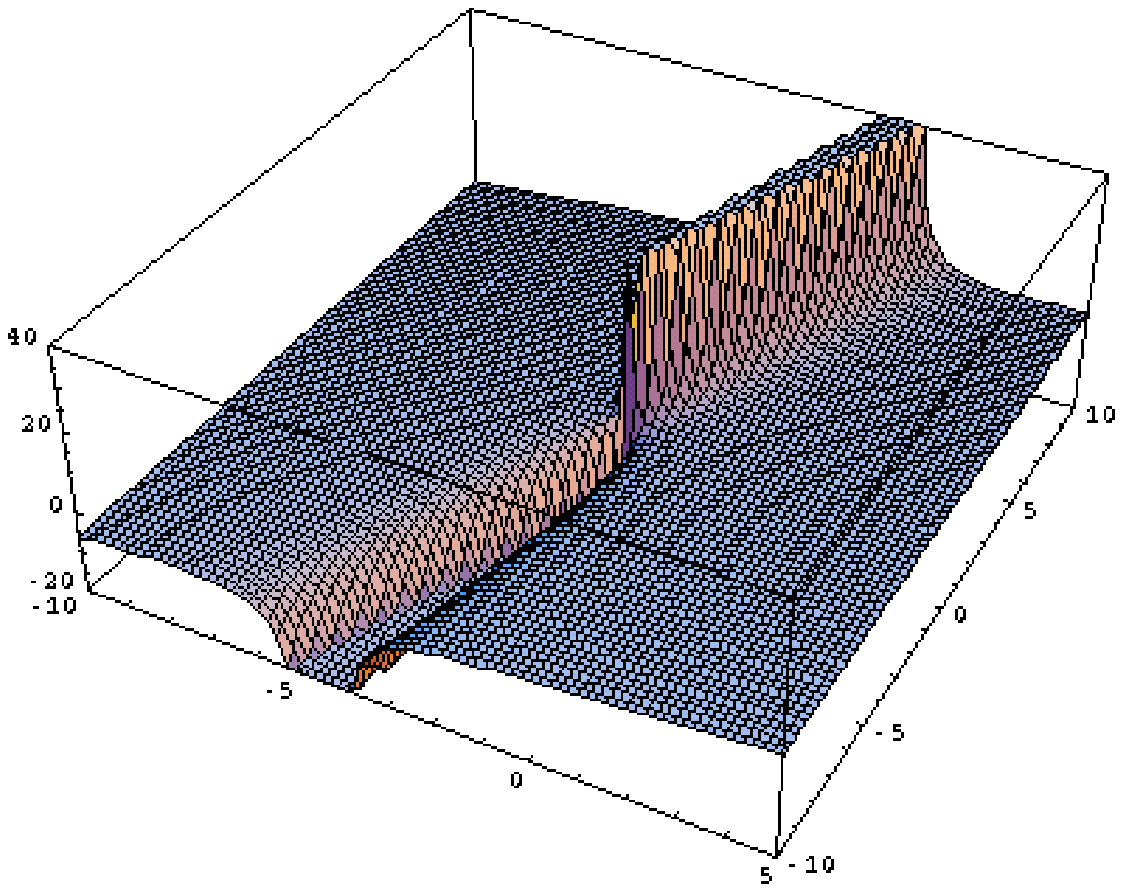}}}

\centerline{\scalebox{1}{\includegraphics{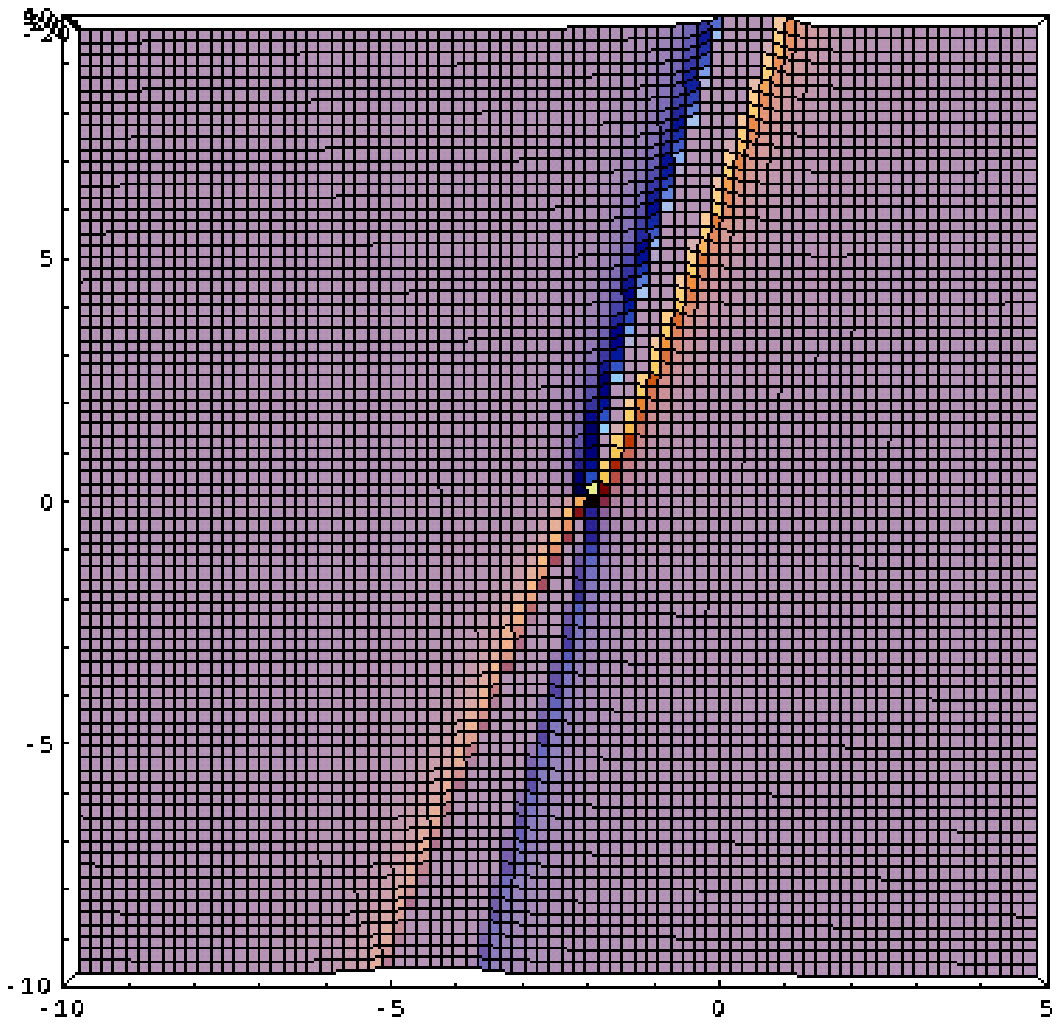}}}

\vspace{2mm}

\centerline{{\bf Figure 2:} $\ds \frac{w_1(x,t)}{\sigma}$ with $N = 1$, $p_1 = 1$, $s_1 = 2$  and 
$c = 0$.} 

\newpage

\centerline{\scalebox{1}{\includegraphics{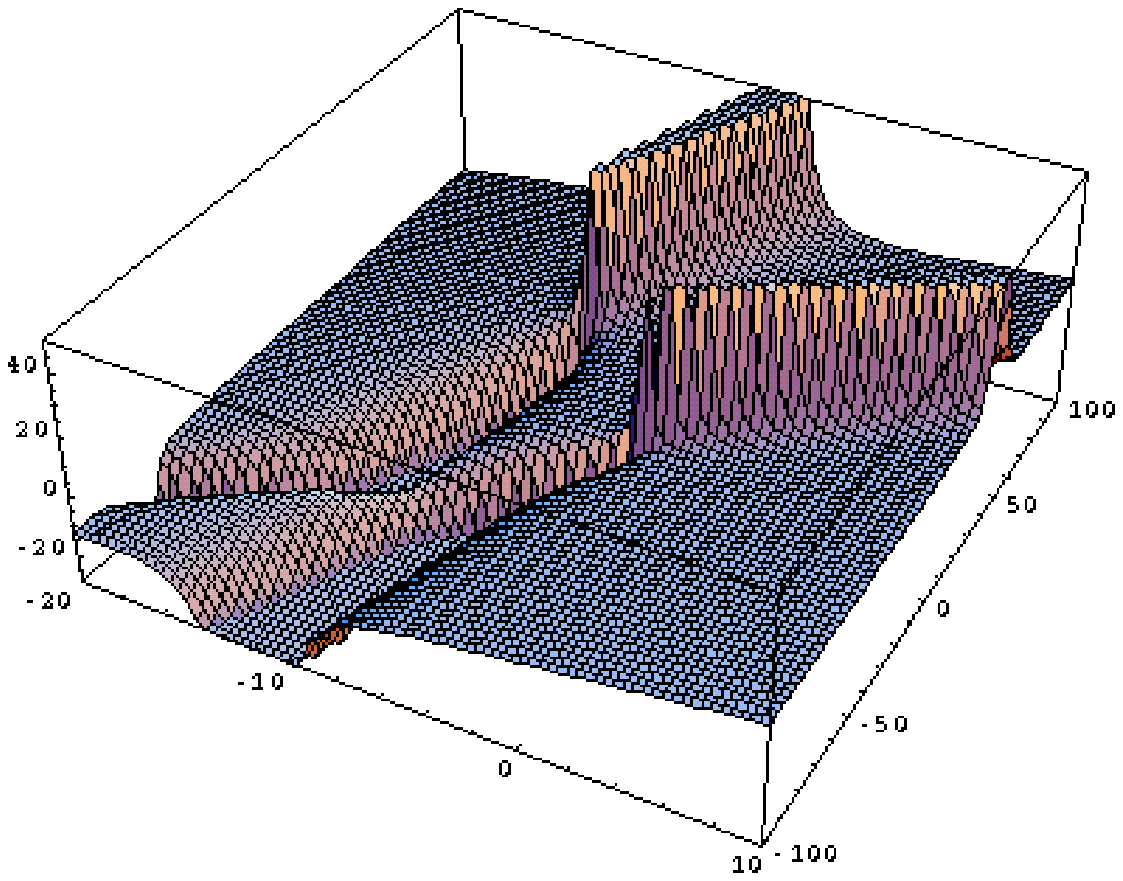}}}

\centerline{\scalebox{1}{\includegraphics{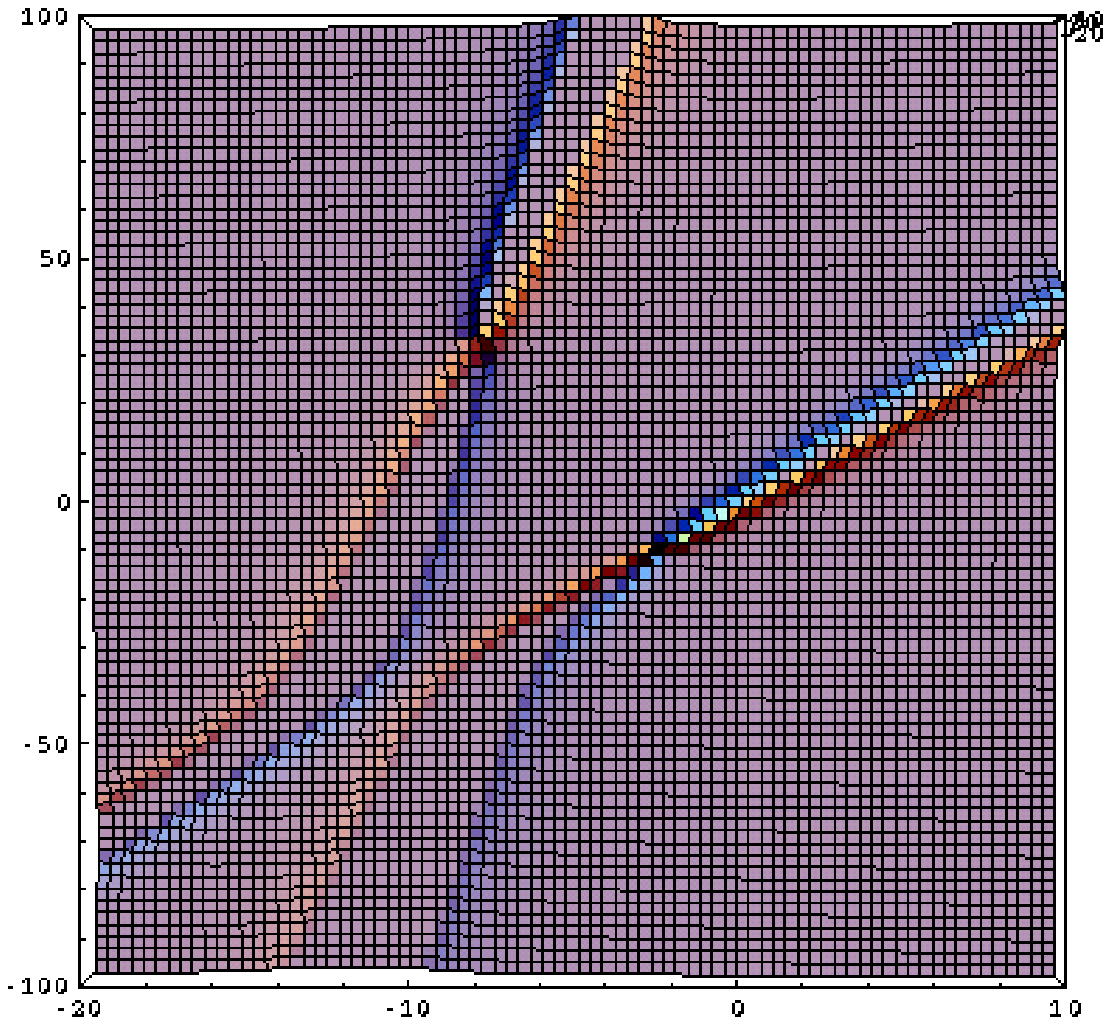}}}

\vspace{2mm}

\centerline{{\bf Figure 3:} $\ds \frac{w_2(x,t)}{\sigma}$ with $N = 2$, $p_1 = 1$, $s_1 = 2$, 
$p_2 = 0.5$, $s_2 = 5$ and $c = 0$.}

\newpage

\centerline{\scalebox{1}{\includegraphics{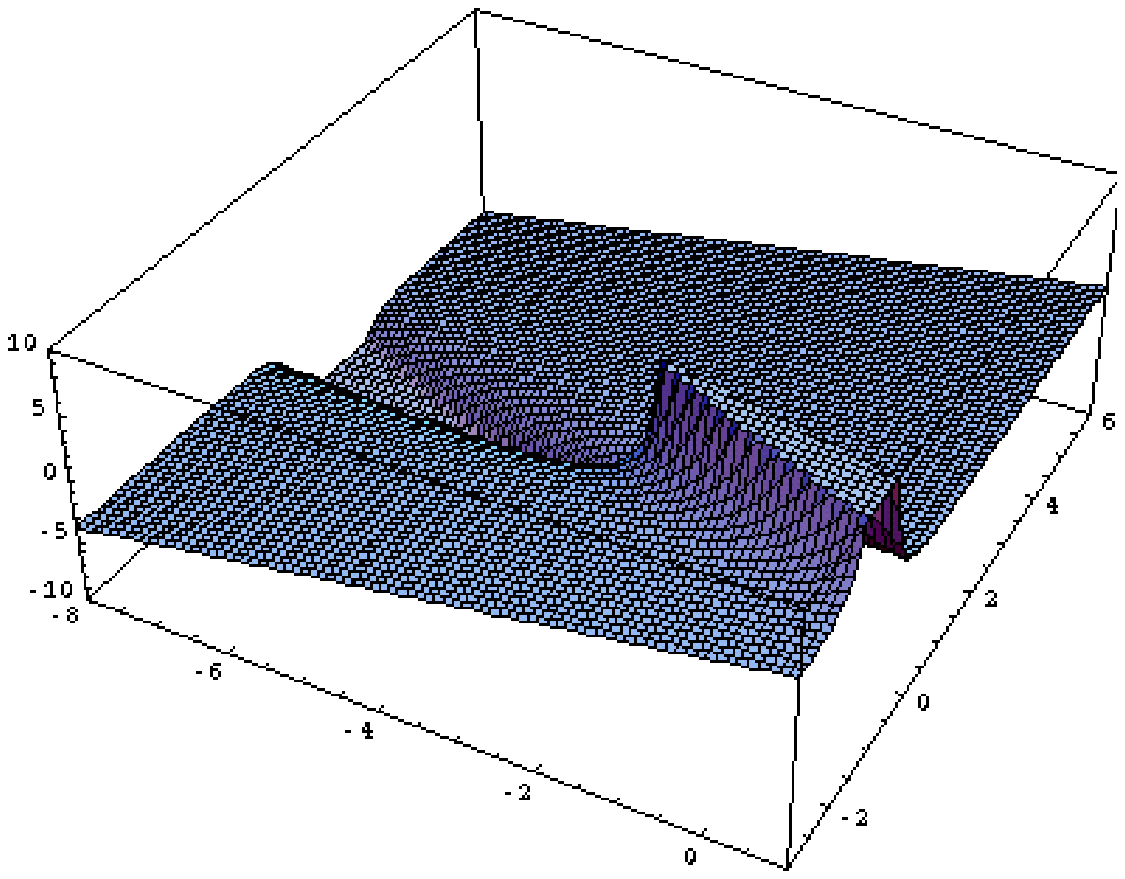}}}

\centerline{\scalebox{1}{\includegraphics{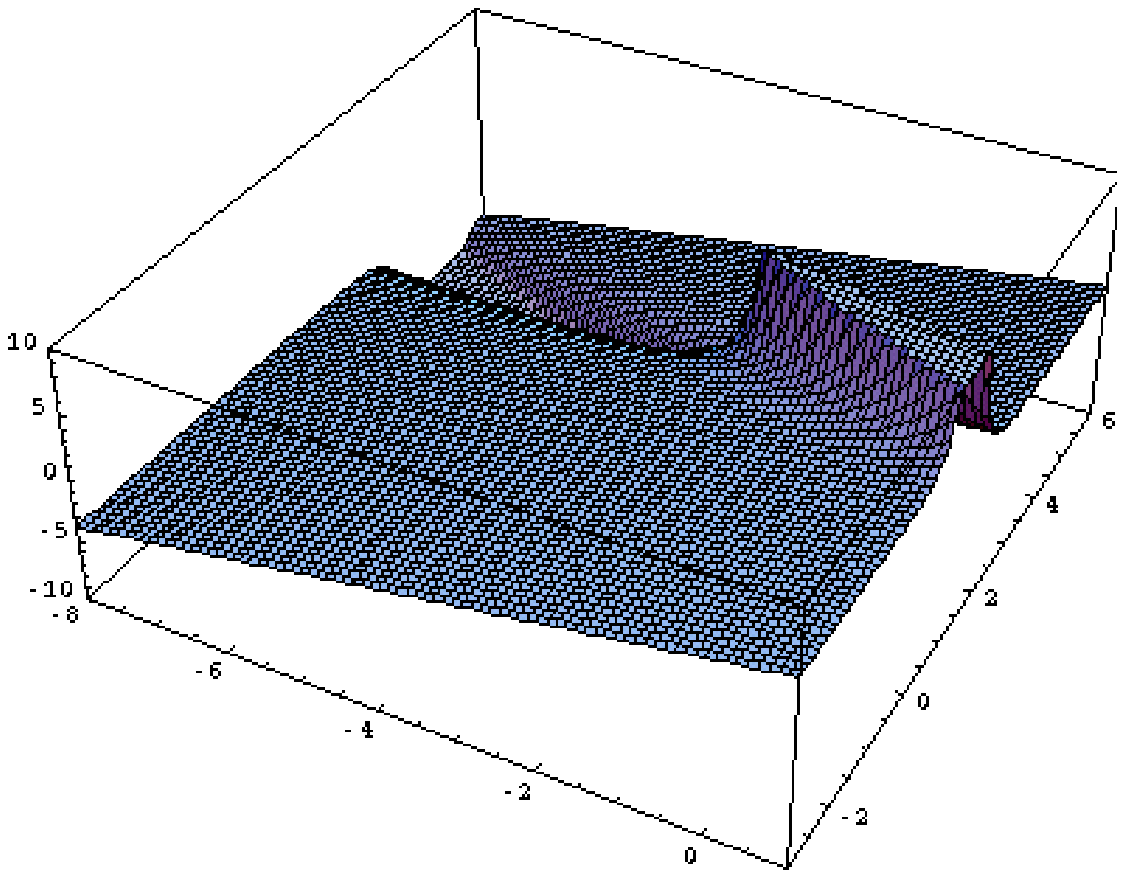}}}

\noindent
{{\bf Figure 4:} Time evolution of $\ds \frac{w_1(x,z,t)}{\sigma}$ with $N = 1$, 
$p_1 = 1$, $q_1 = 3$, $s_1 = 2$ and $c = 0$.}

\newpage

\centerline{\scalebox{1}{\includegraphics{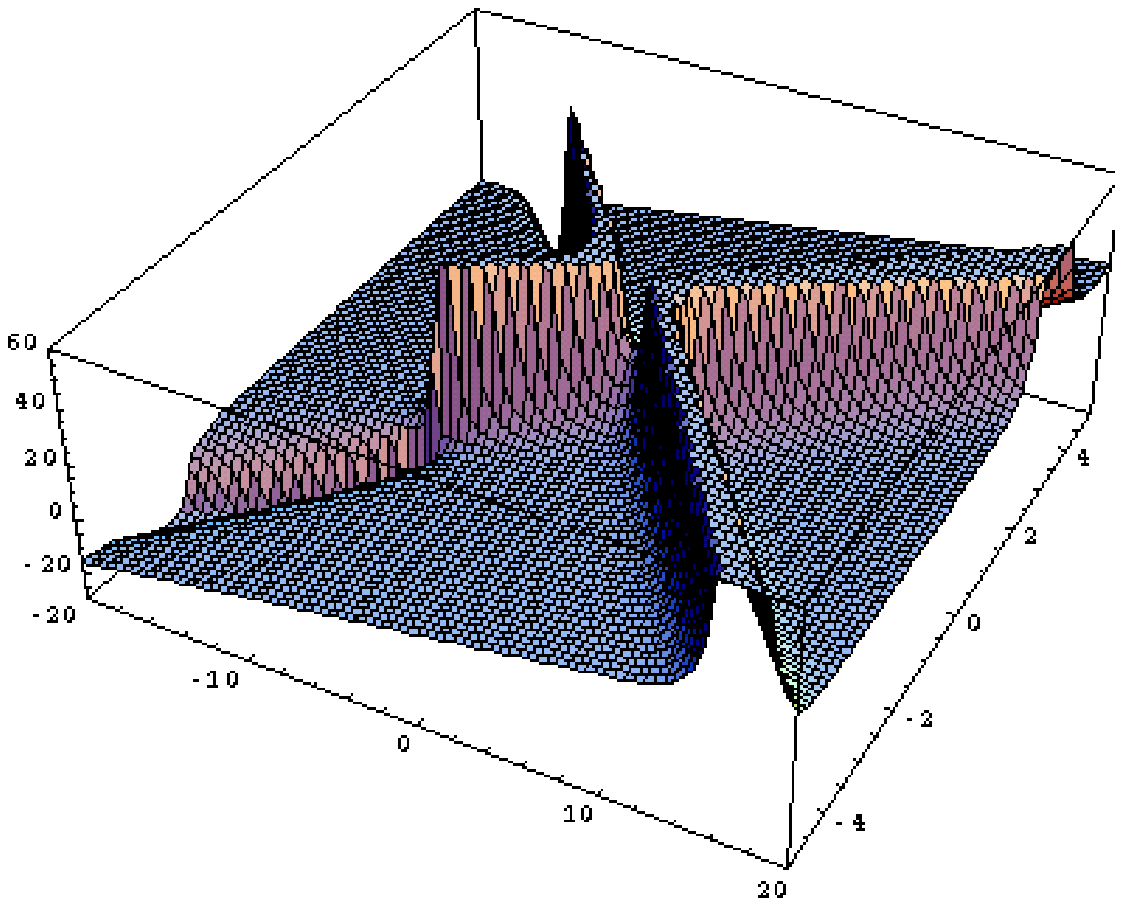}}}

\centerline{\scalebox{1}{\includegraphics{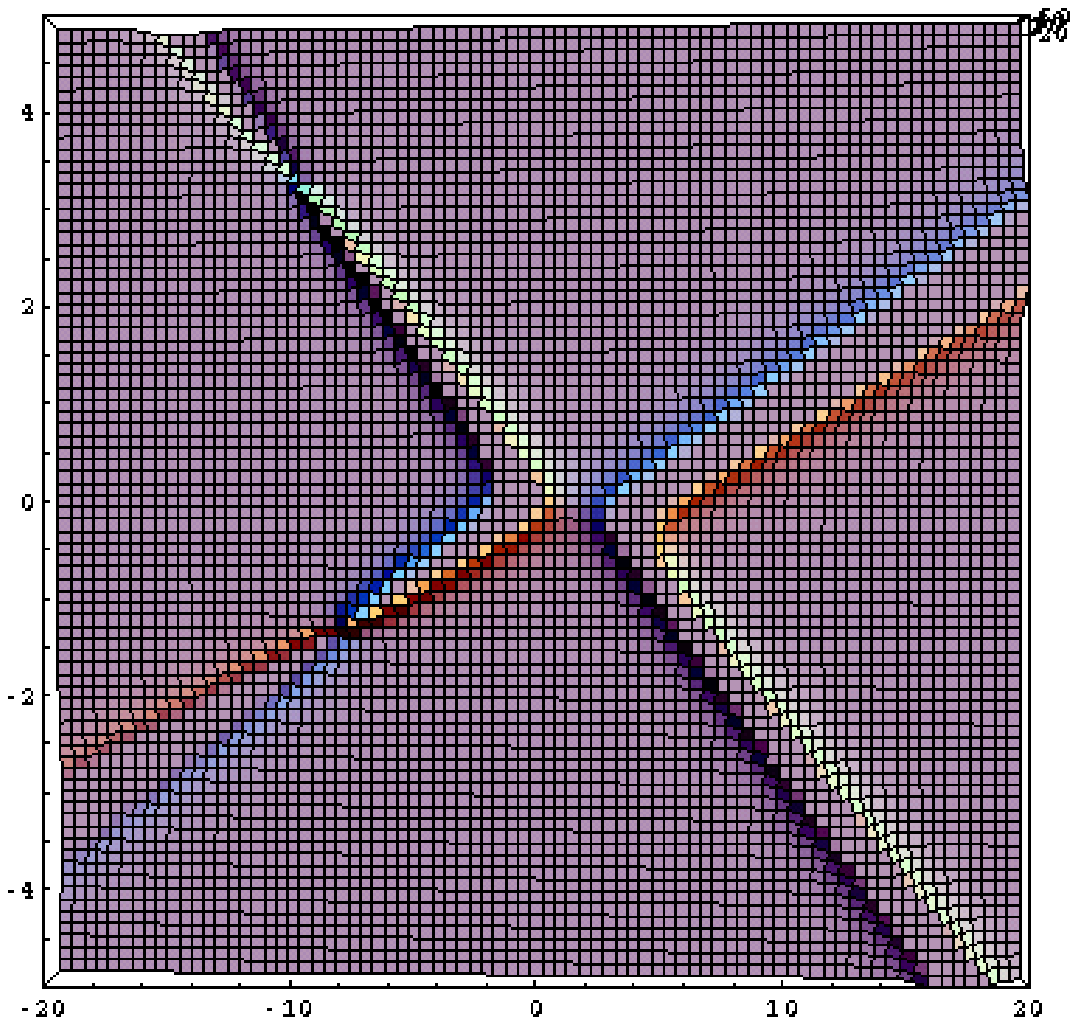}}}

\vspace{2mm}

\noindent
{{\bf Figure 5:} $\ds \frac{w_2(x,z,t)}{\sigma}$ with $N = 2$, $p_1 = 1$, $q_1 = 3$, 
$s_1 = 0$, $p_2 = 0.5$, $q_2 = - 3$, $s_2 = 0$, $c = 0$ and $t = 0$.}

\newpage

\subsection*{Acknowledgements}
A part of this work was done at Yukawa Institute for Theoretical Physics, 
University of Kyoto, Japan. Numerical computations were carried 
out at the Yukawa Institute Computer Facility.  We
would like to thank T.~Fukuyama, R.~Kubo, 
N.~Sasa and R.~Sasaki, for useful discussions. We also thank an 
anonymous referee for valuable remarks.

\label{yu-toda-lp}

\end{document}